\documentclass[12pt]{amsart}
\newtheorem{thm}{Theorem}[section]
\newtheorem{lemma}[thm]{Lemma}
\newtheorem{coro}[thm]{Corollary}
\newtheorem{rem}[thm]{Remark}
\newtheorem{defn}[thm]{Definition}
\newtheorem{prop}[thm]{Proposition}
\newcommand{\R}{{\mathbb R}}

\newcommand{\Z}{{\mathbb Z}}
\newcommand{\C}{{\mathbb C}}
\newcommand{\T}{{\mathbb T}}
\newcommand{\Q}{{\mathbb Q}}
\newcommand{\nn}{\nonumber}
\newcommand{\bea}{\begin{eqnarray}}
\newcommand{\eea}{\end{eqnarray}}
\newcommand{\ba}{\begin{array}}
\newcommand{\ea}{\end{array}}

\newcommand{\ol}{\overline}

\newcommand{\bth}{\begin{thm}}
\newcommand{\eth}{\end{thm}}
\newcommand{\bl}{\begin{lemma}}
\newcommand{\el}{\end{lemma}}
\newcommand{\bk}{\begin{coro}}
\newcommand{\ek}{\end{coro}}
\newcommand{\bh}{\begin{hypo}}
\newcommand{\eh}{\end{hypo}}
\newcommand{\br}{\begin{rem}}
\newcommand{\er}{\end{rem}}
\newcommand{\bpf}{\begin{proof}}
\newcommand{\epf}{\end{proof}}

\begin{document}
\title{Gabor analysis, Noncommutative Tori and Feichtinger's algebra}

\author{Franz Luef}
\address{Fakult\"at f\"ur Mathematik\\
Nordbergstrasse 15\\ 1090 Wien\\ Austria}
\email{Franz.Luef@univie.ac.at}


\maketitle
\pagestyle{myheadings}
\markboth{\small Franz Luef}{\small Gabor Analysis, Noncommutative Tori and Feichtinger's algebra } 

\begin{center}
  \sffamily{\it
"Often in mathematics, understanding comes from generalization,\\
instead of considering the object {\sl per se} when one tries to\\
find the concepts which embody the power of the object."}
\end{center}
\rightline{Alain Connes}
\section{Abstract}
We point out a connection between Gabor analysis and noncommutative analysis.
Especially, the strong Morita equivalence of noncommutative tori appears as
underlying setting for Gabor analysis, since the construction of
equivalence bimodules for noncommutative tori has a natural formulation in
the notions of Gabor analysis. As an application we show that Feichtinger's
algebra is such an equivalence bimodule. Furthermore, we present Connes's
construction of projective modules for noncommutative tori and the relevance
of a generalization of Wiener's lemma for twisted convolution by Gr\"ochenig
and Leinert. Finally we indicate an approach to the biorthogonality relation
of Wexler-Raz on the existence of dual atoms of a Gabor frame operator based
on results about Morita equivalence.

\section{Introduction}

By definition, Morita equivalent algebras $\mathcal{A}$ and $\mathcal{B}$
have equivalent categories of (left) modules. The main theorem of algebraic
Morita theory states that Morita equivalences are implemented by a
$(\mathcal{A},\mathcal{B})$-{\it bimodule} $X$. In \cite{Mor58} Morita proved
the fundamental theorem, showing that this bimodule is  
invertible if and only if the
bimodule is projective and finitely generated as a left $\mathcal{A}$-module
and as a right $\mathcal{B}$-module and $\mathcal{A}\to End_{\mathcal{B}}(X)$
and $\mathcal{B}\to End_{\mathcal{A}}(X)$ are algebra isomorphisms. For
further discussion on the  Morita equivalence of algebras and the relevant modules
and functors we refer to the recent paper \cite{BW04}.
\par
In the seminal papers \cite{Rief74a,Rief74b,Rief76} Rieffel
developed the notion of {\it strong Morita Equivalence}, 
which is an extension of Morita equivalence of algebras to the setting of
$C^*$-algebras. The relevant category of modules over a
$C^*$-algebra, to be preserved under strong Morita equivalence, is
that of Hilbert spaces on which the $C^*$-algebra acts through
bounded operators. More precisely, for a given $C^*$-algebra
$\mathcal{A}$, we consider the category $Herm(\mathcal{A})$ whose
objects are pairs $(\mathcal{H},\rho)$, where $\mathcal{H}$ is a
Hilbert space and $\rho:\mathcal{A}\to\mathcal{B}(\mathcal{H})$ is
a nondegenerate $*$-homomorphism of algebras, and morphisms are
bounded intertwiners. Since we are dealing with more elaborate
modules, it is natural that a bimodule giving rise to a functor
$Herm(\mathcal{B})\to Herm(\mathcal{A})$ should be equipped with
extra structure. If $(\mathcal{H},\rho)\in Herm(\mathcal{B})$ and
$_\mathcal{A}V_\mathcal{B}$ is an
$(\mathcal{A},\mathcal{B})$-bimodule, then if $V$ itself is
equipped with an inner product $\langle.,.\rangle_{\mathcal{B}}$
with values in $\mathcal{B}$. More precisely, let $V$ be a right
$\mathcal{B}$-module. Then a $\mathcal{B}$-valued inner product
$\langle.,.\rangle_{\mathcal{B}}$ on $V$ is a $\C$-sesquilinear
pairing $V\times V\to\mathcal{B}$ (linear in the second argument)
such that, for all $f_1,f_2\in V$ and $T\in \mathcal{B}$, we have
\begin{enumerate}
  \item  $\langle f_1,f_2\rangle_{\mathcal{B}}=\langle
          f_2,f_1\rangle_{\mathcal{B}}^*$
  \item  $\langle f_1,f_2T\rangle_{\mathcal{B}}=\langle f_1,f_2\rangle_{\mathcal{B}}T$
  \item  $\langle f_1,f_1\rangle_{\mathcal{B}}>0$ if $f_1\ne 0$.
\end{enumerate}
Inner products on left modules are defined analogously, but
linearity is required in the first argument. One can show that
$\|f\|_{\mathcal{B}}:=\|\langle f,f\rangle_{\mathcal{B}}\|^{1/2}$
is a norm in $V$. A (right) {\bf Hilbert} $\mathcal{B}$-{\bf
module} is a (right) $\mathcal{B}$-module $V$ together with a
$\mathcal{B}$-valued inner product
$\langle.,.\rangle_{\mathcal{B}}$ so that $V$ is complete with
respect to $\|.\|_{\mathcal{B}}$. Our treatment of strong Morita
equivalence is largely based on a recent paper of Bursztyn and
Weinstein about the connection of Poisson geometry and
noncommutative geometry, \cite{BW04}.
\par
In the study of wavelets Rieffel/Packer and Woods have recently used Hilbert $C^*$-modules, 
see \cite{PR03,PR04,W04}. In the context of Gabor analysis we mention the work of 
Casazza/Lammers \cite{CL03}, but their work is unrelated to our investigations. 

\par
The main goal of this paper is that Rieffel's results on the
strong Morita equivalence of noncommutative tori have a natural
interpretation in terms of Gabor analysis. Recall, that
\begin{equation*}
T_xf(t)=f(t-x)
\end{equation*}
 denotes a {\it translation} by $x\in\R^d$ for $f$ in $L^2(\R^d)$ and
\begin{equation*}
 M_{\omega}f(t)=e^{2\pi i\omega t}f(t)
\end{equation*}
is a {\it modulation} by $\omega\in\widehat{\R}^d$.
More precisely, we show -- following Rieffel --
\cite{Rief88} that the $C^*$-algebra $C^*(\Lambda)$ generated by
time-frequency shifts $\pi(\lambda)=M_{\omega}T_x$ for
$\lambda=(x,\omega)$ in a lattice $\Lambda$ of the time-frequency
plane $\R^d\times\widehat{\R}^d$ is Morita equivalent to the
$C^*$-algebra of time-frequency shifts $C^*(\Lambda^0)$ generated
by $\pi(\lambda^0)$ for $\lambda^0$ in the {\it adjoint lattice}
$\Lambda^0$, see Section 2 for the definition of the adjoint
lattice.
\par
Our main result is that {\it Feichtinger's algebra} $S_0(\R^d)$ is
a bimodule for the $C^*$-algebras $(C^*(\Lambda),C^*(\Lambda^0))$
or an {\it equivalence} bimodule as called by Rieffel.
\par
Rieffel used Schwartz's space $\mathcal{S}(\R^d)$ in his
construction of a bimodule for the $C^*$-algebras
$(C^*(\Lambda),C^*(\Lambda^0))$. Therefore, our investigations are
a further realization of a general strategy of Feichtinger,
that consists in using
$S_0(\R^d)$ is a good substitute for $\mathcal{S}(\R^d)$ 
in many situations. 

As examples we state Feichtinger's seminal paper
\cite{F81}, Reiter's work on Weil's construction of the
metaplectic group \cite{Rei89} and Gr\"ochenig/ Leinert's result on
the irrational case conjecture in Gabor analysis \cite{GL04}.
The present
paper may be seen as a companion to \cite{GL04}, which provides
some reasons for the relevance of noncommutative tori in Gabor
analysis. On the other hand we show that Rieffel obtained the
Janssen representation of Gabor frame operators, and  
the Fundamental Identity of Gabor analysis for functions in the Schwartz-Bruhat
space $\mathcal{S}(G)$ for a locally compact abelian group $G$ and
a closed subgroup $D$ in $G\times\widehat{G}$ already in 1988. More 
concretely, Rieffel observed that
\begin{equation}\label{A-inner product}
  \langle f,g\rangle_{\Lambda}=\sum_{\lambda\in\Lambda}\langle
  f,\pi(\lambda)g\rangle\pi(\lambda)
\end{equation}
is a $C^*(\Lambda)$-valued inner product for $f,g$ in
$\mathcal{S}(\R^d)$.

It comes quite unexpectedly that Rieffel's work on Morita equivalence
for noncommutative tori deals with the same mathematical
structure as Janssen's work on Gabor frame operators, \cite{Jan95}. 
The main goal of this paper is an exploration of this
observation.
\par
The paper is organized as follows: In Section 3 we discuss
operator algebras of time-frequency shifts for a separated set and
we show that there is a correspondence between symmetry of the
separated set in $\R^d\times\widehat{\R}^d$ and the structure of
the corresponding operator algebra. In this context we review some
well-known results on projective representations and we define the
adjoint set of an arbitrary separated set in
$\R^d\times\widehat{\R}^d$ and discuss its relation to the
structure of the commutant of the operator algebra.
\par
In Section 4 we define noncommutative tori and Feichtinger's
algebra. In this context we show that \eqref{A-inner product}
is a Hilbert $C^*(\Lambda)$-valued inner product for $f,g$ in
$S_0(\R^d)$ using functorial properties of $S_0(\R^d)$.
Furthermore we introduce the Short-Time Fourier Transform and its
properties for functions in Feichtinger's algebra $S_0(\R^d)$.
\par
In Section 5 we prove our main theorem, by showing 
that Feichtinger's algebra is a bimodule with respect to the pair of  $C^*$-algebras $C^*(\Lambda)$ and $C^*(\Lambda^0))$.
We discuss the notion of strong Morita equivalence in this setting
and its relation to Gabor expansions, Gabor frame operators,
$\Lambda$-invariant operators for a lattice $\Lambda$ in
$\R^d\times\widehat{\R}^d$. In the proof of Rieffel's
associativity condition we observe its equivalence to the
Fundamental Identity of Gabor Analysis and the existence of
Janssen representation for Gabor frame operators, \cite{Rief88}.
 Our presentation follows \cite{Rief88}, where he used the Poisson's summation formula for
symplectic Fourier transform with respect to a lattice and its adjoint
lattice in the proof of the Fundamental Identiy of Gabor Analysis.
Furthermore, we present the  Gr\"ochenig/Leinert \cite{GL04} 
result on Wiener's lemma for twisted convolution in the context of Connes' 
construction of projective modules over $C^*(\Lambda)$.
\par
In Section 6 we apply our results on strong Morita equivalence of
the $C^*$-algebras $(C^*(\Lambda),C^*(\Lambda^0))$ to the 
structure of Gabor frames. Especially we derive the Wexler-Raz
biorthogonality principle for the existence of dual windows 
for Gabor frames from the relation between traces 
on the Morita equivalent $C^*$-algebras
on $C^*(\Lambda)$ and $C^*(\Lambda^0)$, respectively.

\section{Operator Algebras of Time-Frequency Shifts}

The most general framework in Gabor analysis builds reconstruction
formulas for functions on $\R^d$ from a set of time-frequency
shifts $\mathcal{G}=\{\pi(X_j)~:~X_j\in A\}$ for a countable
subset $A=\{X_j=(x_j,\omega_j)~:~j\in J\}$ with
$\inf_{j,k}|X_j-X_k|>\delta>0$ in the time-frequency plane
$\R^d\times\widehat{\R}^d$. We will refer to this setting as {\it
non-uniform Gabor systems}.
\par
First we  
derive some results about the structure of the operator
algebra of a non-uniform Gabor system $\mathcal{G}$ from general
principles. First the {\it commutant} $\mathcal{G}'$ of all
bounded operators on $L^2(\R^d)$ commuting with every operator in
$\mathcal{A}$ is a unital Banach algebra of time-frequency shifts
with respect to operator composition and the operator norm,
because $\mathcal{A}$ is a subset of unitary (bounded) operators
on $L^2(\R^d)$. Furthermore, if $\mathcal{G}$ is generated by a
set $A$ which is symmetric with respect to the origin, then our
operator algebra $\mathcal{G}$ is invariant under taking adjoints,
moreover $\mathcal{G}$ is actually a $C^*$-algebra of
time-frequency shifts acting on $L^2(\R^d)$. We have the
following chains between $\mathcal{G}$ and its $n$-th commutant
$\mathcal{G}^{(n)}$:

\begin{eqnarray*}\label{com-chain}
   \mathcal{G}&\subset&\mathcal{G}^{''}=\mathcal{G}^{(4)}=\cdots=\\
   &\mathcal{G}^{'}&=\mathcal{G}^{(3)}=\mathcal{G}^{(5)}=\cdots .
\end{eqnarray*}
From now one, let our set  
of time-frequency shifts $A$ be a lattice
$\Lambda$ of $\R^d\times\widehat{\R}^d$ then
$\mathcal{G}=\mathcal{G}^{''}$. The associated Gabor system is
called to be {\it regular}.
These observations 
indicate a strong relationship between the
symmetry of the set $A$ and the structure of the corresponding
operator algebra $\mathcal{G}$ and its commutant
$\mathcal{G}^{'}$. In the following we discuss the connection
between the set $A$ of points in $\R^d\times\widehat{\R}^d$ and
the set $A^0$ of points, which generate the commutant
$\mathcal{G}^{'}$. This investigation yields a beautiful
structure of regular Gabor systems, which is the very reason for
the existence of all of our subsequent results. Therefore, we treat
the relation between $A$ and $A^0$ more concretely.
\par
Our presentation of the relation between $A$ and $A^0$ relies on
the properties of projective representations of the time-frequency
plane $\R^d\times\widehat{\R}^d$. Later we will need some results
about bicharacters and 2-cocycles associated to our projective
representations. Therefore, we now recall their definitions and
some of their basic properties. For a similar treatment of
projective representation, see \cite{DV04}.
\par
Let $G$ be a locally compact abelian group and $\T$ the
multiplicative group of complex numbers of absolute value $1$. A
map $\beta$ of $G\times G$ with values in $\T$
is a {\it multiplier} or (2-){\it cocycle} if it  satisfies
for all $x,y,z \in G$:
\begin{eqnarray*}
  &&\beta(x,0)=\beta(0,x)=1,\\
  &&\beta(x+y,z)\beta(x,y)=\beta(x,y+z)\beta(y,z).
\end{eqnarray*}
Two cocycles $\beta$
and $\beta'$ are called {\it equivalent} or {\it cohomologous} if
there is a Borel map $c$ of $G$ into $\T$, such that for all $x,y \in G$
\begin{equation*}
  \beta'(x,y)=\beta(x,y)\frac{c(x+y)}{c(x)c(y)}
\end{equation*}

A {\it projective representation} $\pi$ of $G$ is a map of $G$
into the unitary group of a Hilbert space $\mathcal{H}$ such that
for a cocycle $\beta$
\begin{equation*}
  \pi(x)\pi(y)=\beta(x,y)\pi(x+y),~~~\pi(0)=1.
\end{equation*}
\par
The map $X=(x,\omega)\mapsto \pi(x,\omega)=M_{\omega}T_x$ for
$X\in\R^d\times\widehat{\R}^d$ into $\mathcal{U}(L^2(\R^d))$ is a
projective representation of the time-frequency plane with cocycle
$\beta'(X,Y)=e^{2\pi iy\cdot\omega}$ for $X=(x,\omega)$ and
$Y=(y,\eta)$.
\par
This projective representation of the time-frequency plane is
equivalent to a projective representation $\pi$ with
$\underline{\text{the}}$ symplectic bicharacter
$\beta(X,Y)=e^{2\pi i(y\omega-x\eta)}$ via the map
$c:(x,\omega)\mapsto e^{2\pi ix\cdot\omega}$ for $X=(x,\omega)$
and $Y=(y,\eta)$ in $\R^d\times\widehat\R^d$.
\par
Recall that a {\it bicharacter} of a locally compact abelian group
$G$ is continuous map $b$ of $G\times G$ into $\T$, which is a
character in each argument. Any such $b$ induces a morphism
$\gamma=\gamma_b$ of $G$ into $\widehat G$ by
\begin{equation*}
  \langle x,\gamma_b(y)\rangle=b(x,y).
\end{equation*}
A bicharacter is called {\it antisymmetric} if it satisfies
for all $x,y \in G$
\begin{equation*}
  b(x,y)b(y,x)=1,~~~b(x,x)=1
\end{equation*}
and {\it symplectic} if it is alternating and $\gamma_b$ is an
isomorphism.
\par
In the study of projective representations $\pi$ of $G$ with
cocyle $\beta$ symplectic bicharacters of $G$ appear naturally in
the commutation rule
\begin{equation}\label{commutation-rel}
  \pi(x)\pi(y)\pi(x)^{-1}\pi(y)^{-1}=b(x,y)I
\end{equation}
with $b(x,y)=\frac{\beta(x,y)}{\beta(y,x)}$ for $x,y\in G$.
\par
For our projective representation by time-frequency shifts of
$\R^d\times\widehat{\R}^d$ we recover the bicharacter
$\rho(X,Y)=e^{2\pi i(y\omega-\eta x)}$ for $X=(x,\omega)$ and
$Y=(y,\eta)$. Due to the close relation of the projective
representation of $\R^d\times\widehat{\R}^d$ by time-frequency
shifts $\pi(x,\omega)=M_{\omega}T_x$ with Heisenberg's commutation
relation the bicharacter $\rho(X,Y)=e^{2\pi i(y\omega-\eta x)}$ is
called the {\it Heisenberg cocycle} for
$\R^d\times\widehat{\R}^d$.
\par
The commutation relation (\ref{commutation-rel}) motivates the
following definition.
\begin{defn}
  Let $A$ be a subset of $G$ and $b$ a bicharacter of $G$. Then the
  {\it adjoint set} $A^0_b$ of $A$ with respect to $b$ is given by
  \begin{equation*}
    A^0_b=\{x\in G~:~b(x,a)=1~~\text{for all}~~a\in A\}.
  \end{equation*}
\end{defn}

The adjoint of a set $A$ is a closed subgroup of $G$ by continuity
of the character $b$. Furthermore a subgroup $A$ of $G$ is called
to be {\it isotropic} for $b$ if $b|_{A\times A}\equiv 1$, or
equivalently $A\subset A^0_b$. We call a group $A$ of $G$ {\it
maximal isotropic} for $b$ if $A=A^0_b$, which are again closed
subgroups of $G$.
\par
We are now in position to answer our original question on the
relation between a set $A$ generating a non-uniform Gabor system
$\mathcal{G}$ and the set $A^0$ generating the commutant
$\mathcal{G}^{'}$.
\begin{eqnarray*}
   A&\subset&A^{00}_b=A^{(4)}_b=\cdots=\\
   &A^{0}_b&=A^{(3)}_b=A^{(5)}_b=\cdots ,
\end{eqnarray*}
where $A^{(n)}_b$ denotes the $n$th adjoint of $A$.
\par
The above chains of relations are the set analogous of
\eqref{com-chain} and therefore, the map $A\mapsto A^{(0)}_b$ is
the desired correspondence.
\par
For the projective representation of $\R^d\times\widehat{\R}^d$ by
time-frequency shifts $\pi(x,\omega)=M_{\omega}T_x$ with the
Heisenberg cocycle $\rho$ we obtain the well-known {\it adjoint
group} $\Lambda^0$ of a regular Gabor system generated by a
lattice $\Lambda$ of $\R^d\times\widehat{\R}^d$. In this case the
maximal isotropic lattice of $\R^d\times\widehat{\R}^d$ is the
standard (von Neumann) lattice $\Z^d\times\Z^d$.
\par
Let us mention that the operator algebra
$\mathcal{G}=\{\pi(\lambda):\lambda\in\Lambda\}$ is a
commutative group of time-frequency shifts if and only if
$\Lambda$ is isotropic, i.e. $\Lambda\subset \Lambda^0$.
\par 
It is an important fact that one can interpret
 the Heisenberg cocycle as follows:  $\rho$ is a character of
$\R^d\times\widehat{\R}^d$, and   every character of
$\R^d\times\widehat{\R}^d$ is of the form
\begin{equation}
  X\mapsto\rho(X,X')~~~\text{for
  some}~~~X'\in\R^d\times\widehat{\R}^d.
\end{equation}
This induces an isomorphism between $\R^d\times\widehat{\R}^d$ and
its dual group $\widehat{\R}^d\times\R^d$.
\par
If $\Lambda$ is a lattice in  
$\R^\times\widehat{\R}^d$,  then every character of $\Lambda$
extends to a character of $\R^d\times\widehat{\R}^d$ and therefore
every character of $\Lambda$ is of the form
\begin{equation}
  \lambda\mapsto\rho(\lambda,Y), \quad ~~~\lambda\in\Lambda
\end{equation}
for some $Y\in\R^d\times\widehat{\R}^d$, where $Y$ needs not to be
unique. The homomorphism from $\R^d\times\widehat{\R}^d$ to
$\widehat{\Lambda}$ has as kernel the {\it adjoint lattice}
\begin{equation}
  \Lambda^0=\{Y\in\R^d\times\widehat{\R}^d~|~~\rho(\lambda,Y)=1
  ~~\text{for all}~~\lambda\in\Lambda\}.
\end{equation}
As a consequence we have, that the adjoint set
of a lattice has the structure of a lattice.
The skew-bicharacter
$\rho$ of $\R^d\times\widehat{\R}^d$ gives a Fourier transform
${\widehat F}^s$ on the time-frequency plane. We call
\begin{eqnarray}\nonumber
{\widehat
 F}^s(Y)&=&\iint_{\R^d\times\widehat{\R}^d}\rho(Y,X)F(X)dX\\\nonumber
                 &=&\iint_{\R^d\times\widehat{\R}^d}e^{2\pi i(y\cdot\omega -
  x\cdot\eta)}F(x,\omega)dxd\omega
\end{eqnarray}
the {\it symplectic Fourier transform} of $F\in
L^2(\R^d\times\widehat{\R}^d)$, since it is induced by the
symplectic form $\Omega$ on $\R^d\times\widehat{\R}^d$. The
symplectic Fourier transform will be essential in our proof of the
Fundamental Identity of Gabor Analysis.

\section{Noncommutative Tori and Feichtinger's Algebra}

A {\it noncommutative} $2d$-{\it torus} $\mathcal{A}_{\Theta}$ is
the universal $C^*$-algebra generated by $2d$ unitaries
$U_1,...,U_{2d}$ subject to the commutation relations
\begin{equation*}
  U_jU_k=e^{2\pi i\theta_{jk}}U_kU_j,\hspace{1.5cm}k,j=1,...,2d
\end{equation*}
for a skew symmetric matrix $\Theta=\big(\theta_{jk}\big)$ with
real entries.
\par
We regard $\Theta$ as a real skew-bilinear form on $\Z^{2d}$, with
entries given by $\Theta(e_j,e_k)=\theta_{jk}$. Then a
noncommutative $2d$-torus $\mathcal{A}_{\Theta}$ is the
twisted group $C^*$-algebra $C^*(\Z,\beta)$, where
$\beta:\Z^{2d}\times\Z^{2d}\to\T$ is a 2-cocycle such that
\begin{equation*}
  \beta(\lambda,\lambda')\overline{\beta(\lambda',\lambda)} = e^{2\pi
i\Theta(\lambda,\lambda')} \quad \quad ~~~\text{for}~~~\lambda,\lambda'\in\Z^{2d} \, .
\end{equation*}
We remark that a noncommutative $2$-torus $A_{\theta}$ is often
called {\it rotation algebra}. The commutation rules for the two
unitary operators $U$ and $V$ generating $A_{\theta}$ read as
\begin{equation*}
  UV=e^{2\pi i\theta}VU \, ,
\end{equation*}
for a real number $\theta$. Let $\theta=\alpha\beta$ then the
$C^*$-algebra generated by time-frequency shifts $\pi(\alpha
k,\beta l)=M_{\beta l}T_{\alpha k}$ for $k,l\in\Z^{d}$ is a
representation of $\mathcal{A}_{\theta}$ on $L^2(\R^d)$.
\par
In this paper we want to treat the general case of
$2d$-noncommutative tori or equivalently $C^*$-algebras
$C^*(\Lambda,\beta)$ of time-frequency shifts $\pi(\lambda)$ for a
lattice $\Lambda$ in $\R^d\times\widehat{\R}^d$ with
\begin{equation*}
  \pi(\lambda)\pi(\lambda')=\beta(\lambda,\lambda')\pi(\lambda+\lambda'), \quad ~~~
  \lambda,\lambda'\in\Lambda.
\end{equation*}

Therefore, an element of $C^*(\Lambda,\beta)$ is given by
\begin{equation*}
  \sum_{\lambda\in\Lambda}a_{\lambda}\pi(\lambda)
\end{equation*}
for an arbitrary
complex-valued sequence ${\bf
a}=(a_{\lambda})_{\lambda\in\Lambda}$.
\par
Moreover, this representation is {\bf faithful} on $L^2(\R^d)$,
which is of great importance in our proofs. One consequence, is
that it is sufficient to establish statements for a dense subspace
of $L^2(\R^d)$. For an operator algebraic proof see \cite{Rief88}
and in \cite{GL01} a proof is given using time-frequency methods.
\par
The choice of a sequence spaces on $\Lambda$ induces on the
noncommutative torus an additional structure. The space
$\mathcal{S}(\Lambda)$ of sequences on
$\Lambda$ with decay faster than the inverse of any polynomial yields a {\it smooth} structure on the algebra of
functions on $C^*(\Lambda,\beta)$, i.e.
\begin{equation*}
  \mathcal{A}_{\Lambda}^{\infty}=\{A\in\mathcal{B}\big(L^2(\R^d)\big)~:~
  A=\sum_{\lambda}a_{\lambda}\pi(\lambda),~~{\bf
  a}=(a_{\lambda})_{\lambda\in\Lambda}\in\mathcal{S}(\Lambda)\}.
\end{equation*}
In the present paper we introduce another structure on
$C^*(\Lambda,\beta)$. Namely,
\begin{equation}
  \mathcal{A}_{\Lambda}^{1}=\{A\in\mathcal{B}\big(L^2(\R^d)\big)~:~
  A=\sum_{\lambda}a_{\lambda}\pi(\lambda),~~{\bf
  a}=(a_{\lambda})_{\lambda\in\Lambda}\in\ell^1(\Lambda)\}.
\end{equation}
The  commutation rules for the  time-frequency shifts  $\pi(\lambda)$
give $C^*(\Lambda,\beta)$ the following structure:
\begin{enumerate}
  \item Let $A_1=\sum_{\lambda\in\Lambda}a_1(\lambda)\pi(\lambda)$
  and $A_2=\sum_{\lambda\in\Lambda}a_2(\lambda)\pi(\lambda)$ for
  ${\bf a_1},{\bf a_2}\in\ell^1(\Lambda)$ then the product of
  $A_1$ and $A_2$ is given by
  \begin{equation*}
    A_1\cdot A_2=\sum_{\lambda\in\Lambda}{\bf a_1}\natural_{\Lambda}{\bf
  a_2}(\lambda)\pi(\lambda),
  \end{equation*}
where
\begin{equation*}
  {\bf a_1}\natural_{\Lambda}{\bf
  a_2}(\lambda)=\sum_{\mu\in\Lambda}a_1(\mu)a_2(\lambda-\mu)\beta(\mu,\lambda-\mu)
\end{equation*}
denotes {\it twisted convolution} of ${\bf a_1}$ and ${\bf a_2}$ and ${\bf a_1}\natural_{\Lambda}{\bf
  a_2}(\lambda)$ is again in $\in\ell^1(\Lambda)$.  
\item Let $A=\sum_{\lambda\in\Lambda}a(\lambda)\pi(\lambda)$ for
${\bf a}\in\ell^1(\Lambda)$ then involution
\begin{equation*}
  A^*=\sum_{\lambda\in\Lambda}a(\lambda)^*\pi(\lambda)
\end{equation*}
induces an involution on $\ell^1(\Lambda)$:
\begin{equation}
  a(\lambda)^*=\beta(\lambda,\lambda)\overline{a(-\lambda)}.
\end{equation}
\end{enumerate}
Therefore, we have that
\begin{prop}
  $\big(\ell^1(\Lambda),\natural_{\Lambda},*\big)$ is an involutive Banach 
   algebra.
\end{prop}
For the construction of a Hilbert $C^*(\Lambda,\beta)$-module $V$ we
are looking for a time-frequency homogenous Banach space,
with properties similar to $\mathcal{S}(\R^d)$. In his seminal paper
\cite{F81} Feichtinger has introduced such a space, nowadays called {\it
Feichtinger's algebra} and denoted by $S_0(\R^d)$. 

Hence, we recall the definition of Feichtinger's algebra and some of
its basic properties. For a detailed discussion
we refer the reader to \cite{F81,FK98,CharlyBook}. There
are many characterizations of Feichtinger's algebra $S_0(\R^d)$. 
The connection between Feichtinger's algebra and Rieffel's work is
most transparent in terms of time-frequency analysis.
\par
Time-frequency representations contain information about the
content of time and frequency in a signal $f$. For our further
investigations we restrict our considerations to the {\it Short
Time Fourier Transform} (STFT). The STFT of a function $f$ with
respect to a window function $g$ is defined for $f,g\in L^2(\R^d)$
as
\begin{equation}
  V_gf(x,\omega)=\int_{\R}f(t)\overline{g(t-x)}e^{-2\pi
  it\cdot\omega}dt,~~~(x,\omega)\in\R^d\times\widehat{\R}^d
\end{equation}
or equivalently as
\begin{equation}
  V_gf(x,\omega)=\langle f,M_{\omega}T_x g\rangle,~~~
  \text{for}~~(x,\omega)\in\R^d\times\widehat{\R}^d.
\end{equation}
Feichtinger's algebra $S_0(\R^d)$ is defined as follows
\begin{equation}
 S_0(\R^d)=\{f\in L^2(\R^d)~~|~~\iint_{\R^d\times\widehat{\R}^d}
 |V_{\varphi}f(x,\omega)|dxd\omega<\infty\},
\end{equation}
where $\varphi(x)=2^{-d/4}e^{-\pi x\cdot x}$ is a Gaussian and its
norm is defined by  
\begin{equation*}
  \|f\|_{S_0}=\iint_{\R^d\times\widehat{\R}^d}|V_{\varphi}f(x,\omega)|dxd\omega.
\end{equation*}
Any other non-zero Schwartz function $g\in\mathcal{S}(\R^d)$, instead of the
Gaussian $\varphi$, yields the same space and an equivalent norm for $S_0(\R)$.
\begin{rem}
Feichtinger's algebra is a particular example of a class of Banach
spaces, the so-called {\it modulation spaces}, which Feichtinger
defined via integrability and decay conditions on the STFT over
$\R^d\times\widehat{\R}^d$ (cf. \cite{FeiMod03}). 
They have been recognized as the
correct class of function spaces for questions in time-frequency
analysis, especially Gabor analysis, \cite{FG89a,FG97}.
\end{rem}
\begin{thm}
  $S_0(\R^d)$ is a Banach algebra under pointwise multiplication.
\end{thm}
By this definition of $S_0(\R)$ elementary properties of the
STFT yield invariance properties of $S_0(\R)$. In the following
lemma we state two well-known facts about STFT.
\begin{lemma}\label{STFTprop}
  Let $f,g\in L^2(\R^d)$ and $(u,\eta)\in\R^d\times\widehat{\R}^d$. Then
  \begin{enumerate}
    \item {\it Covariance Property} of the STFT
          \begin{equation*}\label{covarianceSTFT}
            V_g(\pi(u,\eta)f)g(x,\omega)=e^{2\pi i
           u\cdot(\omega-\eta)}V_gf(x-u,\omega-\eta).
          \end{equation*}
    \item
          \begin{equation*}\label{FI-STFT}
            V_gf(x,\omega)=e^{-2\pi ix\cdot\omega}V_{\hat{g}}\hat{f}(\omega,-x).
          \end{equation*}
  \end{enumerate}
\end{lemma}
\begin{proof}
\indent
  \begin{enumerate}
    \item The covariance property of the STFT is a consequence of the
    of the commutation relation
      \begin{equation*}
        T_xM_{\omega}=e^{-2\pi
        ix\cdot\omega}M_{\omega}T_x,~~~(x,\omega)\in\R^d\times\widehat{\R}^d,
      \end{equation*}
      furthermore one has by definition of the STFT:  
      \begin{eqnarray*}
         V_g(\pi(u,\eta)f)g(x,\omega)&=&\langle M_{\eta}T_u f,M_{\omega}T_xg
         \rangle\\
         &=&\langle f,T_{-u}M_{-\eta}M_{\omega}T_xg\rangle\\
         &=&e^{2\pi i u\cdot(\omega-\eta)}V_gf(x-u,\omega-\eta).
      \end{eqnarray*}
        \item The formula expresses, that for the STFT the Fourier transform yields a rotation
              of the time-frequency plane by an angle of
              $90^{\circ}$. Its another manifestation of the fact
              that STFT contains information of $f$ and $\hat f$.
     \end{enumerate}
\end{proof}
Many properties of Feichtinger's algebra $S_0(\R^d)$ are
elementary consequences of properties of the STFT. The following
theorem may be seen as a realization of this principle, where we
show translation invariance and Fourier invariance of $S_0(\R^d)$
from Lemma (\ref{STFTprop}).
\begin{thm}
  Let $f\in S_0(\R^d)$ and $(u,\eta)\in\R^d\times\widehat{\R}^d$. Then
    \begin{enumerate}
      \item $\pi(u,\eta)f\in S_0(\R^d)$ and
            $\|f\|_{S_0}=\|\pi(u,\eta)f\|_{S_0}.$
      \item $\hat{f}\in S_0(\R^d)$ and $\|\hat{f}\|_{S_0}=\|f\|_{S_0}.$
    \end{enumerate}
\end{thm}
\begin{proof}
  \indent
  \begin{enumerate}
    \item By definition of $S_0(\R^d)$ we have by the invariance of the Gaussian $\varphi$ under Fourier transform that
        \begin{eqnarray*}
          \|f\|_{S_0}&=&\iint_{\R^d\times\widehat{\R}^d}|V_{\varphi}
          (M_{\eta}T_uf)(x,\omega)|dxd\omega\\
          &=&\iint_{\R^d\times\widehat{\R}^d}|V_{\varphi}f(x-u,\omega-\eta)|dxd\omega.
        \end{eqnarray*}
    \item For the invariance under the Fourier transform we use
    (2) of Lemma \eqref{STFTprop}, and that the definition of $S_0(\R^d)$ is
    independent of the window $g\in\mathcal{S}(\R^d)$ and that different windows yield
    equivalent norms for $S_0(\R^d)$:
    \begin{eqnarray*}
      \|\hat f\|_{S_0(\R^d)}&=&\iint_{\R^d\times\widehat{\R}^d}
      |V_{\varphi}\hat f(x,\omega)|dxd\omega\\
      &\le&C\iint_{\R^d\times\widehat{\R}^d}|V_{\hat\varphi}\hat f(x,\omega)|dxd\omega\\
      &=&C\iint_{\R^d\times\widehat{\R}^d}|V_{\varphi}f(-\omega,x)|dxd\omega\\
      &=&C\iint_{\R^d\times\widehat{\R}^d}|V_{\varphi}f(x,\omega)|dxd\omega=\|f\|_{S_0}.
    \end{eqnarray*}
  \end{enumerate}
\end{proof}
Later in our treatment of Rieffel's associativity condition for
$\big(C^*(\Lambda),C^*(\Lambda^0)\big)$ we shall need the following
properties of $S_0(\R^d)$.
\begin{thm}
  Let $f,g$ in $S_0(\R^d)$ then
  \begin{equation*}
    V_gf\in S_0(\R^d\times\widehat{\R}^d).
  \end{equation*}
\end{thm}
For a proof of this statement we refer the reader to \cite{FK98},
but it also follows from the functorial properties and minimality
of $S_0(\R^d)$ in \cite{Fei81}.
\par
Let $F$ be a function on the time-frequency plane
$\R^d\times\widehat{\R}^d$ then the {\it sampling operator} for a
lattice $\Lambda$ in $\R^d\times\widehat{\R}^d$ is defined as
follows
\begin{equation*}
  {\bf R}_{\Lambda}:F\mapsto\big(F(\lambda)\big)_{\lambda\in\Lambda}.
\end{equation*}
\begin{rem} 
We will also write occasionally  $F|_{\Lambda}$ instead of ${\bf R}_{\Lambda}F$ .
\end{rem}
\begin{thm}
  Let $\Lambda$ be a lattice in $\R^d\times\widehat{\R}^d$ and
  $F\in S_0(\R^d\times\widehat{\R}^d)$ then
  \begin{equation}
    {\bf R}_{\Lambda}F\in S_0(\Lambda)=\ell^1(\Lambda)
  \end{equation}
  and ${\bf R}_{\Lambda}$ is bounded on $S_0(\Lambda)$.
\end{thm}
As a consequence of the last theorem we obtain:
\begin{coro}\label{stft-s0}
  Let $f,g\in S_0(\R^d)$ then $V_gf|_{\Lambda}\in\ell^1(\Lambda)$,
  i.e.
  \begin{equation}
    \sum_{\lambda\in\Lambda}|\langle
    f,\pi(\lambda)g\rangle|<\infty.
  \end{equation}
\end{coro}
\begin{rem}
  For $\Lambda$ discrete $S_0(\Lambda)=\ell^1(\Lambda)$.
\end{rem}
We now define a left action of $C^*(\Lambda,\beta)$ on $S_0(\R^d)$
by a {\it Gabor expansion} for the window $g$ and the lattice
$\Lambda\in\R^d\times\widehat{\R}^d$:
\begin{equation*}
  {\bf a}g=\sum_{\lambda\in\Lambda}a(\lambda)\pi(\lambda)g,\hspace{1cm}{\bf
  a}=(a_{\lambda})\in S_0(\Lambda,\beta).
\end{equation*}
The invariance of $S_0(\R^d)$ under time-frequency shifts implies
 that this action
is well-defined on $S_0(\R^d)$.
\begin{prop}\label{bound-time-fre}
  Let ${\mathbf a}\in S_0(\Lambda,\beta)$ and $g\in S_0(\R^d)$, then
  \begin{equation*}
   \big\|\sum_{\lambda\in\Lambda}a(\lambda)\pi(\lambda)g\big\|_{S_0}\le
   \|{\mathbf a}\|_1\|g\|_{S_0}.
  \end{equation*}
\end{prop}
\begin{proof}
\begin{eqnarray*}
  \big\|\sum_{\lambda\in\Lambda}a(\lambda)\pi(\lambda)g\big\|_{S_0}&\le&
  \sum_{\lambda\in\Lambda}|a(\lambda)|\|\pi(\lambda)g\|_{S_0}\\
  &=&\sum_{\lambda\in\Lambda}|a(\lambda)|\|g\|_{S_0}\\
  &=&\|{\mathbf a}\|_1\|g\|_{S_0}.
\end{eqnarray*}
\end{proof}
\par \noindent
Corollary \ref{stft-s0} and Proposition \ref{bound-time-fre}
yield that for $f,g\in S_0(\R^d)$ the left action of
\begin{equation*}
\langle f,g\rangle_{\Lambda}=V_gf|_{\Lambda}
\end{equation*}
is well-defined on $S_0(\R^d)$. In Gabor analysis the mapping of
$\langle f,g\rangle_{\Lambda}\mapsto\langle f,g\rangle_{\Lambda}g$
is called a Gabor type frame operator with  window $g$ and lattice
$\Lambda$,  denoted by
\begin{equation*}
 S_{g,\Lambda}f=\sum_{\lambda\in\Lambda}\langle
 f,\pi(\lambda)g\rangle\pi(\lambda)g.
\end{equation*}
The above discussion shows that the Gabor type frame operator
$S_{g,\Lambda}$ is continuous on $S_0(\R^d)$ for $f,g\in
S_0(\R^d)$. In Section 6 we present some consequences for the
reconstruction of square-integrable functions $f \in L^2(\R^d)$
and is
\par
Rieffel's central observation was that
\begin{equation}\label{Lambda-inner}
  \langle f,g\rangle_{\mathcal{A}}=\sum_{\lambda\in\Lambda}\langle
  f,\pi(\lambda)g\rangle\pi(\lambda),\hspace{0.5cm}f,g\in\mathcal{S}(\R^d)
\end{equation}
is a $C^*$-valued innerproduct for
$\mathcal{A}=C^*(\Lambda,\beta)$. In the sequel we prove that
\eqref{Lambda-inner} defines a $C^*$-valued innerproduct for
$C^*(\Lambda,\beta)$ for $f,g$ in Feichtinger's algebra
$S_0(\R^d)$.
\par
First we show that \eqref{Lambda-inner} is compatible with the
action of $S_0(\Lambda,\beta)$ on $S_0(\R^d)$. More precisely, we
prove the following proposition.
\begin{prop}
Let $f,g\in S_0(\R^d)$ and let ${\bf a}\in S_0(\Lambda,\beta)$.
Then
  \begin{equation*}
    \langle{\bf a} f,g\rangle_{\Lambda}={\bf a}\natural_{\Lambda}
          \langle f,g\rangle_{\Lambda}.
  \end{equation*}
\end{prop}


\begin{proof} 
  \item For $\lambda\in\Lambda$ we have
\bea\nn
  \langle{\bf a}f,\pi(\lambda)g\rangle&=&\sum_{\lambda'\in\Lambda}{a}(\lambda')\langle\pi(\lambda')f,\pi(\lambda)g
   \rangle\\\nn
   &=&\sum_{\lambda'\in\Lambda}{a}(\lambda')\langle f,\pi^*(\lambda')
   \pi(\lambda)g\rangle\\\nn
   &=&\sum_{\lambda'\in\Lambda}{a}(\lambda')\langle f,\pi(\lambda-\lambda')g
   \rangle\beta(\lambda',\lambda-\lambda')\\\nn
   &=&{\bf a}\natural_{\Lambda}\langle f,\pi(\lambda)g\rangle,
\eea since
$\pi(\lambda')\pi(\lambda-\lambda')=\beta(\lambda',\lambda-\lambda')\pi(\lambda)$.
\end{proof}
We now  come to the statement of one of our main theorems.
\begin{thm}\label{Hilbert-S0}
  Feichtinger's algebra $S_0(\R^d)$ is a left Hilbert
  $C^*(\Lambda,\beta)$-module with respect to the inner product, given for 
  $f,g$ in $S_0(\R^d)$ by
  \begin{equation*}
    \langle
    f,g\rangle_{\mathcal{A}}=\sum_{\lambda\in\Lambda}V_gf(\lambda)\pi(\lambda) \, .
  \end{equation*}
\end{thm}
The proof of Theorem \ref{Hilbert-S0} is postponed now and will be given after the 
discussion of the {\bf F}undamental {\bf I}dentity of {\bf G}abor
{\bf A}nalysis, because the positivity of the innerproduct is a
direct consequence of FIGA. In Section \ref{Sect-Bimodule} we
derive FIGA from an identity for products of STFT by an
application of Poisson's summation formula for the symplectic
Fourier transform, which requires the adjoint lattice $\Lambda^0$
of the lattice $\Lambda$ in $\R^d\times\widehat{\R}^d$. In Section
\ref{Sect-Bimodule} we study the structure of $C^*(\Lambda)$ and
define a $C^*(\Lambda^0)$-valued innerproduct. The proof of
Rieffel's associativity condition is an elementary reformulation
of the FIGA. Therefore, Section \ref{Sect-Bimodule} is the natural
place for our presentation of the FIGA.
\par
At the end of this section we present generalizations of some
notions of Hilbert spaces to Hilbert $C^*$-modules.
\par
Let $V$ be a Hilbert $\mathcal{A}$-module then a {\it Hilbert
module map} from $V$ to $V$ is a linear map $T:V\to V$ that
respects the module action: $T({\bf a}f)={\bf a}~(T(f))$ 
for ${\bf
a}\in\mathcal{A}$ and $f\in V$.  The adjoint of an operator on a
Hilbert space plays a central role in the study of operators on
Hilbert spaces and of operator algebras, such as $C^*$-algebras or
von Neumann algebras of operators. The following definition gives
a generalization of adjoints to Hilbert $\mathcal{A}$-module.
\begin{defn}
  Let $V$ be a Hilbert $\mathcal{A}$-module. A map $T:V\to V$ is
  {\it adjointable} if there exists a map $T^*:V\to V$ satisfying
  \begin{equation*}
    \langle Tf,g\rangle_{\mathcal A}=\langle f,T^*g\rangle_{\mathcal A}
  \end{equation*}
  for all $f,g$ in $\mathcal{A}$. The map $T^*$ is called the {\it
  adjoint} of $T$. We denote the set of all adjointable maps by
  $\mathfrak{L}(V)$ and the set of all bounded module maps in $V$
  by $\mathfrak{B}(V)$.
\end{defn}
An elementary consequence of the definitions is the following
facts about adjointable maps.
\begin{enumerate}
  \item Let $T$ be in $\mathfrak{L}(V)$ , then its adjoint is unique and
        adjointable with $T^{**}=T$.
  \item Let $T,S$ be in $\mathfrak{L}(V)$, then
        $ST\in\mathfrak{L}(V)$ with $(ST)^*=T^*S^*$.
  \item $\mathfrak{L}(V)$ equipped with the operator norm
        $\|T\|=\sup\{\|Tx\|~:~\|x\|\le 1\}$ is a $C^*$-algebra.
 \item $\mathfrak{B}(V)$ equipped with the operator norm
        $\|T\|=\sup\{\|Tx\|~:~\|x\|\le 1\}$ is a Banach algebra.
\end{enumerate}
In the case of $S_0(\R^d)$ as $C^*(\Lambda,\beta)$-module the
adjointable maps are those operators $T:S_0(\R^d)\to S_0(\R^d)$
where $T^*$ commutes with all time-frequency shifts
$\{\pi(\lambda)~:~\lambda\in\Lambda\}$. By definition of the
$C^*(\Lambda,\beta)$-innerproduct we have
\begin{eqnarray*}
  \langle Tf,g\rangle_{\mathcal{A}}&=&\sum_{\lambda\in\Lambda}\langle
  Tf,\pi(\lambda)g\rangle\pi(\lambda)\\
  &=&\sum_{\lambda\in\Lambda}\langle f,T^*\pi(\lambda)g\rangle\pi(\lambda)\\
  &=&\sum_{\lambda\in\Lambda}\langle f,\pi(\lambda) T^*g\rangle\pi(\lambda)\\
  &=&\langle f,Tg\rangle_{\mathcal{A}}.
\end{eqnarray*}

In \cite{FK98} Feichtinger and Kozek treated selfadjoint operators
on $S_0(\R^d)$, which commute with
$\{\pi(\lambda)~:~\lambda\in\Lambda\}$. They called those
operators $\Lambda$-invariant. The set of all selfadjoint
adjointable operators of $\mathfrak{L}(S_0(\R^d))$ is an ideal in
$\mathfrak{L}(S_0(\R^d))$.
\par
The notion of finite rank operators and of compact module
operators is of great relevance in the construction of Morita
equivalences between $C^*$-algebras, see Section
\ref{Sect-Bimodule}.
\begin{defn}
  Let $f,g$ be elements of a Hilbert $\mathcal{A}$-module $V$.
  Then a {\bf rank one operator} $K_{f,g}:V\to V$ is defined by
  \begin{equation*}
    K_{f,g}h:=\langle f,h\rangle_{\mathcal{A}}g.
  \end{equation*}
 The set of {\bf compact Hilbert module operators} on $V$ is the
 closed subspace of $\mathfrak{L}(V)$ generated by the rank one
 maps $K_{f,g}$. We denote the set of compact Hilbert module operators
 by ${\mathcal K}(V)=\overline{\{K_{f,g}:f,g\in V\}}$.
\end{defn}
\begin{rem}
  A compact Hilbert module operator is not necessarily a bounded
  operator on $V$, but for Hilbert $\C$-modules $\mathcal{H}$ the notion
  specializes to the definition of a compact operator on
  $\mathcal{H}$.
\end{rem}
The following proposition gives some elementary facts about
compact Hilbert module operators.
\begin{prop}
  Let $f,g\in V$ and $T\in\mathcal{L}(V)$. Then one has 
  \begin{enumerate}
    \item $K_{f,g}$ is adjointable and $K_{f,g}^*=K_{g,f}$.
    \item $TK_{f,g}=K_{Tf,g}$.
    \item $K_{f,g}T=K_{f,T^*g}$.
    \item $\|K_{f,g}\|\le\|f\|\|g\|$.
  \end{enumerate}
\end{prop}
A direct consequence of the previous observations is the following
statement.
\begin{prop}
  Let $V$ be a Hilbert $\mathcal{A}$-module,
  $\mathfrak{K}(V)$ is a closed ideal in $\mathfrak{L}(V)$.
\end{prop}
Now we investigate the set of rank one module operators for our
Hilbert $C^*(\Lambda,\beta)$ module $S_0(\R^d)$. By definition a
rank one module operator is given by
\begin{eqnarray*}
  K_{g,f}\gamma&=&\langle f,\gamma\rangle_{\mathcal{A}}g  \\
               &=&\sum_{\lambda\in\Lambda}\langle
               f,\pi(\lambda)\gamma\rangle\pi(\lambda)g=S_{g,\gamma,\Lambda}f,
\end{eqnarray*}
for $f,g,\gamma\in S_0(\R^d)$. The operator $S_{g,\gamma,\Lambda}$
is called a Gabor frame operator with analysis window $\gamma$ and
synthesis window $g$ for a lattice $\Lambda$.

Therefore, a
finite rank module operator is a finite sum of Gabor frame
operators, a so-called {\it multi-window} Gabor frame operator.
Furthermore, a rank one module operator $S_{g,\gamma,\Lambda}$ is
an adjointable operator, i.e. it is a $\Lambda$-invariant
operator. This elementary fact has far reaching consequences, see
Section \ref{gab-wexler-raz}.

\section{Feichtinger's Algebra as Bimodule for $C^*(\Lambda)$ and
$C^*(\Lambda^0)$}\label{Sect-Bimodule}

In Section 3 we discussed the relation between an operator algebra
of time-frequency shifts generated by a lattice $\Lambda$ in
$\R^d\times\widehat{\R}^d$ and its commutant. In this section we
continue discussion in the light of Morita equivalence of
$C^*$-algebras.
\par
The adjoint lattice of $\Lambda$ in $\R^d\times\widehat{\R}^d$ was
defined as the set of all points $X=(x,\omega)$ in
$\R^d\times\widehat{\R}^d$ such that $\rho(\lambda,X)=1$, which by
the commutation relation of time-frequency shifts
\eqref{commutation-rel} is equivalent to
\begin{equation*}
  \Lambda^0=\{\lambda^0\in\R^d\times\widehat{\R}^d:\pi(\lambda)\pi(\lambda^0)
  =\pi(\lambda^0)\pi(\lambda)~~~\text{for
  all}~~\lambda\in\Lambda\}.
\end{equation*}
Therefore, the set of all bounded operators on $L^2(\R^d)$
commuting with elements from $C^*(\Lambda,\beta)$ is the
$C^*$-algebra generated by time-frequency shifts $\pi(\lambda^0)$
for $\lambda^0$ in $\Lambda^0$.
\par
In Section 4 we have defined a left action of $C^*(\Lambda,\beta)$
on $S_0(\R^d)$. Now $S_0(\R^d)$ has the structure of a
bimodule, where the right action is induced by the opposite
algebra of $C^*(\Lambda^0,\beta)$. Following Rieffel in
\cite{Rief88}, $C^*(\Lambda^0,\beta)$ can be generated by
$\pi^*(\lambda^0)$ acting on the left on $S_0(\R^d)$, which
commutes with the right action of $C^*(\Lambda,\beta)$ on
$S_0(\R^d)$. Therefore, the opposite algebra of
$C^*(\Lambda,\beta)^{\text{opp}}$ is generated by
$\pi^*(\lambda^0)$ with
$\overline{\beta}(X,Y)=\overline{\beta(X,Y)}$ for $X=(x,\omega)$
and $Y=(y,\eta)$ as cocycle, i.e. $C^*(\Lambda^0,\ol{\beta})$. By
definition the opposite algebra of $C^*(\Lambda,\beta)$ gives a
right action on $S_0(\R^d)$ by a Gabor expansion with respect to
the lattice $\Lambda^0$
\begin{equation*}
  g{\bf
  b}=|\Lambda|^{-1}\sum_{\lambda^0\in\Lambda^0}b(\lambda^0)\pi^*(\lambda^0)g,~~~
  g\in S_0(\R^d),{\bf b}\in S_0(\Lambda^0,\ol{\beta}).
\end{equation*}
Note, that cohomologous cocycles yield isomorphic $C^*$-algebras.
By a reasoning similar to the one used  in Section 4 for the left action
$C^*(\Lambda,\beta)$ we obtain that the right action is well-defined
on $S_0(\R^d)$.
\par
Before $S_0(\R^d)$ is given the structure of a right
$C^*(\Lambda^0,\ol{\beta})$ we state the Fundamental Identity of
Gabor analysis, because it is essential in our construction of the
bimodule $S_0(\R^d)$ for $C^*(\Lambda,\beta)$ and
$C^*(\Lambda^0,\ol{\beta})$.
\par
\begin{thm}[FIGA]
Let $f_1,g_1,f_2,g_2\in S_0(\R^d)$, then
\begin{equation*}
  \sum_{\lambda\in\Lambda}V_{g_1}f_1(\lambda)\overline{V_{g_2}f_2(\lambda)}=
  |\Lambda|^{-1}\sum_{\lambda^0\in\Lambda^0}V_{g_1}g_2(\lambda^0)
  \overline{V_{f_1}f_2(\lambda^0)}
\end{equation*}
\end{thm}
In \cite{Rief88} Rieffel proved FIGA for Schwartz functions
$f_1,f_2,g_1,g_2$ in $\mathcal{S}(G)$ for an elementary locally
compact abelian group $G$. In \cite{TO} Tolmieri and Orr proved a
special case of Rieffel's result for functions on $\R$ in their
study of Gabor frames. Later, Janssen continued the work of
Tolmieri/Orr and introduced a representation of Gabor frame
operators, Janssen's representation \cite{Jan95}. In his proof of
the Morita equivalence of $C^*(\Lambda,\beta)$ and
$C^*(\Lambda^0,\ol{\beta})$ Rieffel had derived Janssen's
representation of a Gabor frame operator.
\par
Following Rieffel we use Poisson summation formula for symplectic
Fourier transform in the proof of FIGA. The following theorem
states the Poisson summation formula for the symplectic Fourier
transform, see Section 2 for the definition.
\begin{thm}
Let $F\in S_0(\R^d\times{\widehat{\R}}^d)$ then
\begin{equation}
  \sum_{\lambda\in\Lambda}F(\lambda)=|\Lambda|^{-1}\sum_{\lambda^0\in\Lambda^0}
  {\widehat{F}}^s(\lambda^0)
\end{equation}
holds pointwise and with absolute convergence of both sums.
\end{thm}

\begin{proof}(FIGA)\label{FIGA}
\par
 If $f,g\in S_0(\R^d)$ we have that $V_gf\in
 S_0(\R^d\times\widehat{\R}^d)$. Then
 $F=V_{g_1}f_1\ol{V_{g_2}f_2}$ is in
 $S_0(\R^d\times\widehat{\R}^d)$, because
 $S_0(\R^d\times\widehat{\R}^d)$ is a Banach algebra under
 multiplication. Poisson's summation formula for $F$ yields FIGA.
 Therefore, we compute the symplectic Fourier transform of $F$.
 \begin{eqnarray}\nn
   \widehat{F}^s(Y)&=&\iint_{\R^d\times{\widehat{\R}}^d}
                      V_{g_1}f_1(X)\overline{V_{g_2}f_2(X)}\rho(Y,X)dX\\\nn
                   &=&\iint_{\R^d\times{\widehat{\R}}^d}
                   \langle\pi(Y)f_1,\pi(Y)\pi(X)g_1\rangle
                   \overline{\langle
                   f_2,\pi(X)g_2\rangle}\rho(X,Y)dX\\\nn
                   &=&\iint_{\R^d\times{\widehat{\R}}^d}
                   \langle\pi(Y)f_1,\pi(X)\pi(Y)g_1\rangle
                   \overline{\langle
                   f_2,\pi(X)g_2\rangle}\rho(X,Y)dX\\\nn
                   &=&\langle f_1,\pi(Y)f_2\rangle\overline{\langle
                   g_1,\pi(Y)g_2\rangle},
 \end{eqnarray}
 where in the last step we used Moyal's formula.
\end{proof}
As a first application of FIGA we finish the proof of Theorem
\ref{Hilbert-S0} by showing the positivity of
\begin{equation}
  \langle
  f,f\rangle_{\mathcal{A}}=\sum_{\lambda\in\Lambda}
  \langle f,\pi(\lambda)f\rangle\pi(\lambda)
\end{equation}
as an operator on $L^2(\R^d)$.
\begin{prop}
  Let $f\in S_0(\R^d)$ then $\langle f,f\rangle_{\mathcal{A}}$ is
  a positive element of $C^*(\Lambda,\beta)$.
\end{prop}
\begin{proof}
The representation of time-frequency shifts of
$C^*(\Lambda,\beta)$ is faithful, therefore, it suffices to
establish positivity for a dense subspace of $L^2(\R^d)$. Of
course we choose $S_0(\R^d)$ as dense subspace. Let $g\in
S_0(\R^d)$
  \begin{eqnarray}\nn
     \langle\langle f,f\rangle_{\Lambda}g,g\rangle&=&
     \Big\langle\sum_{\lambda\in\Lambda}\langle f,\pi(\lambda)f\rangle\pi(\lambda)g,g\Big\rangle\\\nn
      &=&\sum_{\lambda\in\Lambda}\langle f,\pi(\lambda)f\rangle
         \ol{\langle g,\pi(\lambda)g\rangle}\\\nn
      &=&\sum_{\lambda^0\in\Lambda^0}\langle f,\pi(\lambda^0)g\rangle
         \ol{\langle f,\pi(\lambda^0)g\rangle}\ge 0.
  \end{eqnarray}
\end{proof}
In an analogous manner as in our discussion of Theorem
\ref{Hilbert-S0} we get that the right action of
$C^*(\Lambda^0,\ol{\beta})$ with properly defined
$\natural_{\Lambda^0}$ and involution $*$ defines a right Hilbert
$C^*(\Lambda^0,\ol{\beta})$-module structure on $S_0(\R^d)$ with
respect to the
$\mathcal{B}:=C^*(\Lambda^0,\ol{\beta})$-innerproduct
\begin{equation*}
  \langle f,g\rangle_{\mathcal{B}}:=|\Lambda|^{-1}\sum_{\lambda^0\in\Lambda^0}
  \langle \pi(\lambda^0)g,f\rangle\pi(\lambda^0),\hspace{0.8cm}f,g\in
  S_0(\R^d).
\end{equation*}
Two $C^*$-module structures $\big(\mathcal{A},\langle
.,.\rangle_{\mathcal{A}}\big)$ and $\big(\mathcal{B},\langle
.,.\rangle_{\mathcal{B}}\big)$ on a bimodule $V$ are compatible if
\begin{equation}\label{rieffel-assoc}
  \langle f,g\rangle_{\mathcal{A}}h=f\langle
  g,h\rangle_{\mathcal{B}},\hspace{0.8cm}\text{for all}~~~f,g,h\in V.
\end{equation}
Some authors call \eqref{rieffel-assoc} {\bf Rieffel's
associativity condition} for $\langle.,.\rangle_\mathcal{A}$ and
$\langle.,.\rangle_\mathcal{B}$.
\par
In our setting Rieffel's associativity condition expresses
Janssen's representation of a Gabor frame operator $S_{g,\gamma}$
for a window $g,\gamma\in S_0(\R^d)$.
\begin{thm}
  Let $\mathcal{A}=C^*(\Lambda,\beta)$ and
  $\mathcal{B}=C^*(\Lambda^0,\ol{\beta})$ with the above defined
  actions and innerproducts $\langle.,.\rangle_{\mathcal{A}}$ and
  $\langle.,.\rangle_{\mathcal{B}}$, respectively. Then
  \begin{equation*}
    S_{g,\gamma,\Lambda}f=|\Lambda|^{-1}S_{f,\gamma,\Lambda^0}g
  \end{equation*}
  for all $f,g,\gamma\in S_0(\R^d)$.
\end{thm}
\begin{proof}
  As in the proof of positivity of $\langle
  f,f\rangle_{\mathcal{A}}$ for $f\in S_0(\R^d)$ it suffices to
  show that for all $\gamma, h\in S_0(\R^d)$
  \begin{equation*}
    \langle S_{g,\gamma,\Lambda}f,h\rangle=
    |\Lambda|^{-1}\langle S_{f,\gamma,\Lambda^0}g,h\rangle \, .
  \end{equation*}

  \par
\begin{eqnarray*}
  \big\langle\langle
  f,g\rangle_{\mathcal{A}}\gamma,h\big\rangle&=&\sum_{\lambda\in\Lambda}
  \langle f,\pi(\lambda)g\rangle\ol{\langle h,\pi(\lambda)\gamma\rangle}\\
  &\stackrel{FIGA}{=}&|\Lambda|^{-1}\sum_{\lambda^0\in\Lambda^0}
  \langle f,\pi(\lambda^0)h\rangle\ol{\langle g,\pi(\lambda^0)\gamma\rangle}\\
  &=&\big\langle f\langle
  g,\gamma\rangle_{\mathcal{B}},h\big\rangle.
\end{eqnarray*}
\end{proof}
\par
A Hilbert $C^*$-module $V$ over $\mathcal{A}$ is called {\it full}
when the collection $\{\langle f,g\rangle_{\mathcal{A}}:~f,g\in
V\}$ is dense in $V$.
\begin{defn}
Two $C^*$-algebras $\mathcal{A}$ and $\mathcal{B}$ are {\bf
strongly Morita equivalent} if there exists a full Hilbert
$C^*$-module $V$ over $\mathcal{B}$ such that $\mathcal{B}\simeq
{\mathcal K}(V,\mathcal{A})$.
\end{defn}
\begin{rem}
  We denote by ${\mathcal K}(V,\mathcal{A})$ the $C^*$-algebra of
  compact Hilbert $\mathcal{A}$ operators $K_{f,g}^{\mathcal{A}}$.
\end{rem}
The Morita equivalence of $C^*(\Lambda,\beta)$ and
$C^*(\Lambda^0,\ol{\beta})$ is a consequence of the following
theorem.
\begin{thm}
  Let $S_0(\R^d)$ be given a bimodule structure as defined above. Let $\mathcal{A}=C^*(\Lambda,\beta)$ and
  $\mathcal{B}=C^*(\Lambda^0,\ol{\beta})$ then
  \begin{enumerate}
   \item  $\{\langle f,g\rangle_{\mathcal{A}}:~f,g\in
           S_0(\R^d)\}$ is dense in $\mathcal{A}$, i.e. $S_0(\R^d)$ is a full Hilbert
           $\mathcal{A}$-module.
   \item  $\{\langle f,g\rangle_{\mathcal{B}}:~f,g\in
           S_0(\R^d)\}$ is dense in $\mathcal{B}$, i.e. $S_0(\R^d)$ is a full Hilbert
           $\mathcal{B}$-module.
   \item  For all $f\in S_0(\R^d)$ and $A\in \mathcal{A}$, we have
          \begin{equation*}
            \langle fA ,fA\rangle_{\mathcal{A}}\le\|A\|^2 \langle
            f,f\rangle_{\mathcal{A}},
          \end{equation*}
          i.e. boundedness of the right action.
  \item  For all $f\in S_0(\R^d)$ and $B\in \mathcal{B}$, we have
          \begin{equation*}
            \langle Bf,Bf\rangle_{\mathcal{B}}\le\|B\|^2 \langle
            f,f\rangle_{\mathcal{B}},
          \end{equation*}
          i.e. boundedness of the left action.
  \end{enumerate}
  implies that $S_0(\R^d)$ is an equivalence bimodule
  $(\mathcal{A},\mathcal{B})$ with norm $\|f\|:=\langle f,f\rangle_{\mathcal{A}}^{1/2}$.
\end{thm}
\begin{proof}
Our proof follows Rieffel's approach, see \cite{Rief88}.
\begin{enumerate}
  \item  The linear span of the range of
      $\langle .,.\rangle_{\mathcal{A}}$ is an ideal in $\mathcal{A}$. Then
      the norm closure $I$ of this linear span is an ideal in
      $\mathcal{A}$. Furthermore $I$ is invariant under modulation
      and because $\pi(\lambda)$ is a faithful representation of
      $\mathcal{A}$, we get the desired conclusion.
  \item By similar arguments as in $(1)$.
  \item It suffices to verify the inequality for a dense subspace of $L^2(\R^d)$. Let
  $h\in S_0(\R^d)$ and $A\in C^*(\Lambda,\beta)$, then
  \begin{eqnarray*}
    \big\langle h\langle A
    f,Af\rangle_{\mathcal{A}},h\big\rangle&=&\big\langle\langle h,
    Af\rangle_{\mathcal{B}}A f,h\big\rangle\\
    &=&\langle A f,\langle A f,h\rangle_{\mathcal{B}}h\rangle\\
    &=&\big\langle A f, Af\langle
    h,h\rangle_{\mathcal{A}}\big\rangle\\
    &=&\big\langle A(f\langle h,h\rangle_{\mathcal{A}})^{1/2},
    A(f\langle h,h\rangle_{\mathcal{A}})^{1/2} \big\rangle\\
    &\le&\|A\|^2\big\langle f,f\langle
    h,h\rangle_{\mathcal{A}}\big\rangle\\
    &=&\|A\|^2\big\langle h\langle
    f,f\rangle_{\mathcal{A}},h\big\rangle
  \end{eqnarray*}
  holds for all $f$ in $S_0(\R^d)$. A standard density argument
  yields the desired result.
  \item By similar arguments as in $(3)$.
\end{enumerate}

\end{proof}
\begin{coro}
 $C^*(\Lambda,\beta)$ and $C^*(\Lambda^0,\ol{\beta})$ are {\it strongly Morita equivalent}.
\end{coro}
By definition a {\it projective module} $V$ is isomorphic to a
direct summand of a free module ${\mathcal A}^n$ with standard
basis $\{e_j\}$, i.e. there is a self-adjoint $n\times n$-matrix
$P$ with entries in ${\mathcal A}$ which is a projection, such
that $V=P{\mathcal A}^n$. Rieffel proved that if ${\mathcal A}$
and $\mathcal{B}$ are unital $C^*$-algebras and if $V$ is a
$({\mathcal B},{\mathcal A})$-equivalence bimodule, then $V$ is a
projective right ${\mathcal B}$-module, and a projective left
${\mathcal A}$-module. Furthermore, ${\mathcal A}$ is equivalent
to the $C^*$-algebra ${\mathcal K}(V,{\mathcal B})$ of compact
Hilbert ${\mathcal B}$-module operators.
\par
In particular, let ${\mathcal B}=C^*(\Lambda^0,\ol{\beta})$ and
let $V$ denote the right ${\mathcal A}$-module obtained by
completing $S_0(\R^d)$ as described earlier. Then, we have:
\begin{thm}
  Feichtinger's algebra $S_0(\R^d)$ is a {\bf finitely generated projective} ${\mathcal
  B}$-module and ${\mathcal K}(S_0(\R^d),{\mathcal B})$ is
  equivalent to $C^*(\Lambda,\beta)$.
\end{thm}
In \cite{Rief81a} Rieffel made the observation that finitely generated
projective $C^*$-modules possess a reconstruction formula in terms
of a tight module frame, which is a generalization of the familiar
notion of tight frames for Hilbert spaces. In a subsequent paper
we discuss the connection between tight module frames for
$C^*(\Lambda,\beta)$ and the characterization of $S_0(\R^d)$ with
multi-window Gabor frames.
\par
We finish this section with a connection between Connes'
construction of projective ${\mathcal A}$-modules in \cite{Con81}
and the recent result of Gr\"ochenig and Leinert on Wiener's lemma
for twisted convolution, \cite{GL04}.
\par
Let $\mathcal{A}_{0}$ be the subalgebra of absolutely convergent
series of ${\mathcal A}$, i.e.
\begin{equation*}
   {\mathcal A}_{\Lambda}^1=\{A\in L^2(\R^d):A= {\sum_{{\lambda\in\Lambda}}}
   a(\lambda)\pi(\lambda),~~{\bf a}\in S_0(\Lambda)\}
\end{equation*}
with norm $\|A\|_{{\mathcal A}_{\Lambda}}^1=\|{\bf a}\|_1$. Then the
main result of Gr\"ochenig/Leinert is the following noncommutative
version of Wiener's lemma.
\begin{thm}\label{GroLei}
  If $A\in{\mathcal A}_0$ is invertible in $C^*(\Lambda,\beta)$,
  then $A^{-1}\in\mathcal{A}_0$.
\end{thm}
\par
We present an important result by  Connes,
which due to the complexity of our assumptions is a little bit lengthy.

\begin{thm}[Connes]
Let ${\mathcal A}$ be a unital $C^*$-algebra, let $V$ be a
projective right ${\mathcal A}$-module with ${\mathcal A}$-valued
inner product, and let $B={\mathcal K}({\mathcal A},V)$. Then
${\mathcal B}$ is a $C^*$-algebra and $V$ has a corresponding
${\mathcal B}$-valued innerproduct. Let ${\mathcal A}_0$ and
${\mathcal B}_0$ be dense $*$-subalgebras of ${\mathcal A}$ and
$\mathcal{B}$ respectively containing the identity elements, and
let $V_0$ be a dense subspace of $V$ which is closed under the
actions of ${\mathcal A}_0$ and ${\mathcal B}_0$, and such that
the restrictions to $V_0$ of the innerproduct have values in
${\mathcal A}_0$ and ${\mathcal B}_0$, respectively. If ${\mathcal
B}_0$ has the property that any element of ${\mathcal B}_0$ which
is invertible in ${\mathcal B}$ has its inverse in ${\mathcal
B}_0$, then $V_0$ is a finitely generated projective right
${\mathcal A}_0$-module.
\end{thm}
\par
We now apply this construction to the above given equivalence
bimodule $S_0(\R^d)$ and the lattice $\Lambda$ in$
\R^d\times\widehat{\R}^d$ and $A=C^*(\Lambda^0,\overline{\beta})$
and ${\mathcal B}_0=\{B\in
L^2(\R^d):~~B=\sum_{\lambda^0\in\Lambda^0}
   b(\lambda^0)\pi(\lambda^0),~~{\bf b}\in S_0(\Lambda^0)\}$ then
\begin{thm}\label{fingen}
 $S_0(\R^d)$ is a finitely generated projective right $\mathcal{A}_0$-module.
\end{thm}
For the Schwartz space $\mathcal{S}(\R^d)$ Connes called this a
beautiful result and of great importance in noncommutative
analysis.

\section{Application to Gabor Analysis: Biorthogonality Relation of Wexler-Raz}\label{gab-wexler-raz}

Recently, Gabor frames have been applied in various fields of
mathematics, electrical engineering and signal analysis, see
\cite{FS98,FS03}. In this section we give a first glimpse of the
usefulness of Rieffel's work on strong Morita equivalence of $C^*$-algebras
generated by time-frequency shifts.
\par
Let $\Lambda$ be a lattice in $\R^d\times\widehat{\R}^d$ and $g\in
L^2(\R^d)$ then a {\it Gabor system} ${\mathcal
G}(g,\Lambda):=\{\pi(\lambda)g:~\lambda\in\Lambda\}$ for a Gabor
atom $g\in L^2(\R^d)$ is a {\it Gabor frame} if there are finite
positive reals $A,B$ such that
\begin{equation*}
  A\|f\|^2\le\sum_{\lambda\in\Lambda}|\langle
  f,\pi(\lambda)g\rangle|^2\le B\|f\|^2,~~~\text{for all}~~f\in
  L^2(\R^d).
\end{equation*}
This is equivalent to invertibility and boundedness of the {\it
Gabor frame operator}
\begin{equation*}
  S_{g,\Lambda}f=\sum_{\lambda\in\Lambda}\langle
  f,\pi(\lambda)g\rangle\pi(\lambda)g,~~~\text{for all}~~f\in L^2(\R^d).
\end{equation*}
As a consequence of the invertibility of $S_{g,\Lambda}$ we have
the following reconstruction formulas for $f\in L^2(\R^d)$
\begin{equation}\label{framereconst1}
  f=(S_{g,\Lambda})^{-1}S_{g,\Lambda}f=\sum_{\lambda\in\Lambda}\langle
  f,\pi(\lambda)g\rangle\pi(\lambda)(S_{g,\Lambda})^{-1}g
\end{equation}
or
\begin{equation*}\label{framereconst2}
  f=S_{g,\Lambda}(S_{g,\Lambda})^{-1}f=\sum_{\lambda\in\Lambda}\langle
  f,\pi(\lambda)(S_{g,\Lambda})^{-1}g\rangle\pi(\lambda)g.
\end{equation*}
The coefficients in reconstruction formulas \eqref{framereconst1}
and \eqref{framereconst2} are not unique, because in general
time-frequency shifts $\pi(\lambda)$ and $\pi(\lambda')$ are not
linearly independent for $\lambda,\lambda'\in\Lambda$. Therefore,
many researchers have investigated the set of all possible
dual windows $\gamma$ such that $S_{g,\gamma}=I$. Of special
importance is the function $\gamma_0:=(S_{g,\Lambda})^{-1}g$, the
{\it canonical dual window}. There are many characterizations of
$\gamma_0$ in the set of all dual windows.
\par
The Gabor frame operator $S_{g,\Lambda}$ commutes with
time-frequency shifts $\{\pi(\lambda):\lambda\in\Lambda\}$,
therefore, the dual Gabor frame
$\{\pi(\lambda)\gamma_0:\lambda\in\Lambda\}$ has the structure of
a Gabor frame. This observation and (\ref{framereconst1}) for
$(S_{g,\Lambda})^{-1}f$ yields that the inverse frame operator of
a frame ${\mathcal G}(\Lambda,g)$ is given by
\begin{equation}
  (S_{g,\Lambda})^{-1}f=S_{\gamma_0,\Lambda}f=
  \sum_{\lambda\in\Lambda}\langle
  f,\pi(\lambda)\gamma_0\rangle\pi(\lambda)\gamma_0.
\end{equation}
Gr\"ochenig and Leinert were motivated by a practical question on
the quality of the canonical dual window $S_{g,\Lambda}^{-1}g$ of
a Gabor frame ${\mathcal G}(g,\Lambda)$ generated by a Gabor atom
$g$ in $S_0(\R^d)$. They established that Feichtinger's algebra is
a good class of Gabor atoms. Namely,
\begin{thm}
  Let ${\mathcal G}(g,\Lambda)$ be a Gabor frame generated by $g\in
  S_0(\R^d)$ then the canonical dual window $\gamma_0=S_{g,\Lambda}^{-1}g$
  is in $S_0(\R^d)$.
\end{thm}
\begin{proof}
  Let $g\in S_0(\R^d)$ then by assumption $S_{g,\Lambda}$ is
  invertible in $C^*(\Lambda,\beta)$, but Theorem \ref{GroLei}
  implies that $S_{g,\Lambda}^{-1}$ is an element of ${\mathcal
  A}_{\Lambda}^{1}$, i.e. $\gamma_0\in S_0(\R^d)$.
\end{proof}
In \cite{Jan95} Janssen had proved that for a Gabor frame
${\mathcal G}(g,\alpha\Z\times\beta\Z)$ generated by a Schwartz
function $g$ the canonical dual window is also a Schwartz function
under the restriction that $\alpha,\beta\in\Q$. Janssen had
conjectured that his result is also valid for irrational lattice
constants $\alpha,\beta$. We remark that a resolution of Janssen's
conjecture is a corollary of Connes result that ${\mathcal
S}(\R^d)$ is closed under holomorphic functional calculus,
\cite{Con81} and \cite{Rief88}.
\par
In \cite{WR90} Wexler/Raz characterized the set of all dual atoms with the
structure of  a Gabor frame for Gabor expansions on finite abelian
groups. Their work had been extended to the continuous setting
independently by Daubechies, H.L. Landau and Z. Landau in
\cite{DLL95}, by Janssen in \cite{Jan95} and by Ron and Shen
\cite{RonS93},\cite{RonS97}. In the work on this problem the so-called 
{\it Janssen representation} of a Gabor frame operator was introduced in
\cite{Jan95}. Also Feichtinger and Zimmermann considered this
topic and found the minimal assumptions for the validity of
Wexler-Raz's biorthogonality relation and Janssen's
representation, \cite{FZ98}. In \cite{FK98} and \cite{FZ98}
Feichtinger and his collaborators introduced the notion of the
adjoint lattice for elementary locally compact abelian groups,
which Rieffel already used in his construction of equivalence
bimodules between noncommutative tori, \cite{Rief88}. In this
section we derive the result of Wexler-Raz from the Morita
equivalence of $\C^*(\Lambda,\beta)$ and $C^*(\Lambda^0,\ol{\beta})$ and
the relation between the canonical traces $\tau_{\mathcal{A}}$ and
$\tau_{\mathcal{B}}$, respectively.
\par
One of the early successes of operator algebras was the
classification of all commutative $C^*$-algebras by
Gelfand as the involutive complex-valued continuous functions over
a compact space. Riesz's representation theorem for positive
linear functionals of involutive complex-valued continuous
functions over a compact space $X$ yields to an extension of the
Lebesgue integral. Therefore, integration of continuous functions
over a compact space is considered as a trace on a commutative
$C^*$-algebra. Therefore, traces on general $C^*$-algebras are the
natural framework for non-commutative Radon measure theory.
\par
The existence of canonical traces on ${\mathcal
A}=C^*(\Lambda,\beta)$ and ${\mathcal B}=C^*(\Lambda^0,\ol{\beta})$ is
one of the pleasant properties of noncommutative tori.
\par
First we recall that a {\it trace} $\tau_{\mathcal C}$ on a
$C^*$-algebra ${\mathcal C}$ is a linear functional satisfying
\begin{eqnarray*}
    &&\tau(I)=1,~~\text{for the identity operator}~~I~~\text{of}~~ {\mathcal C},\\
    &&\tau(A B)=\tau(BA),~~\text{for all}~~A,B\in{\mathcal C},\\
    &&\tau(A^*A)>0~~~\text{for all nonzero}~~ A~~\text{in}~~{\mathcal C}.
\end{eqnarray*}
In the case of ${\mathcal A}$ a normalized trace $\tau_{{\mathcal
A}}$ is given by
\begin{equation*}
  \tau_{\mathcal A}(\langle f,g\rangle_{{\mathcal
  A}})=\langle f,g\rangle,\hspace{0.8cm}f,g\in S_0(\R^d),
\end{equation*}
and for ${\mathcal B}$ the canonical trace $\tau_{{\mathcal B}}$
is normalized by
\begin{equation*}
  \tau_{\mathcal B}(\langle f,g\rangle_{{\mathcal
  B}})=|\Lambda|^{-1}\langle f,g\rangle,\hspace{0.8cm}f,g\in S_0(\R^d),
\end{equation*}
which follows from Morita equivalence of ${\mathcal A}$ and
${\mathcal B}$. This fact can be considered as a noncommutative
Poisson summation formula
\begin{equation}\label{noncommPoisson}
  \tau_{\mathcal{A}}(\langle f,g\rangle_{{\mathcal A}})=
  |\Lambda|^{-1}\tau_{{\mathcal B}}(\langle f,g\rangle_{{\mathcal B}}).
\end{equation}

Our restriction in the following theorem to $g,\gamma\in S_0(\R)$
is just for convenience. We refer to Gr\"ochenig's excellent
survey \cite{CharlyBook} of Gabor analysis for the general case of
$g,\gamma\in L^2(\R^d)$.
\begin{thm}[Wexler-Raz]\label{wexler-raz}
The following conditions are equivalent:
  \begin{enumerate}
    \item $S_{g,\gamma}=I$ is a normalized tight frame for $L^2(\R^d)$.
    \item $|\Lambda|^{-1}\langle\gamma,\pi(\lambda^0)g\rangle=\delta_{\lambda,0}$.
  \end{enumerate}
\end{thm}
\begin{proof}
$(2)\Rightarrow(1):$
\par
  Follows from the fact that the identity of
  ${\mathcal A}$ is  $I=\delta_{\lambda,0}$, where
$\delta_{i,k}$ is the Kronecker delta. Therefore, by assumption
\begin{equation*}
\tau_{{\mathcal B}}(S_{g,\gamma})=\tau_{{\mathcal
B}}(I)=\delta_{\lambda,0}
\end{equation*}
and by application of \eqref{noncommPoisson} we get
\begin{equation*}
\tau_{{\mathcal B}}(I)=|\Lambda|^{-1}\tau_{{\mathcal A}}(\langle
g,\gamma\rangle)=|\Lambda|^{-1}\langle\gamma,g\rangle \, .
\end{equation*}
The implication $(1)\Rightarrow(2)$ is trivial in the light of
Rieffel's associativity condition.
\end{proof}
\begin{coro} 
  For $g,\gamma\in S_0(\R^d)$ dual functions the two Gabor systems
  ${\mathcal{G}}(g,\Lambda^0)$ and ${\mathcal{G}}(\gamma,\Lambda^0)$
  are biorthogonal to each other on $L^2(\R^d)$.
\end{coro}
The proof is an elementary reformulation of Theorem
\ref{wexler-raz}.
\begin{coro}
  A Gabor system ${\mathcal{G}}(g,\Lambda)$ is a tight frame if
  and only if ${\mathcal{G}}(\gamma,\Lambda^0)$ is an orthonormal
  system with frame bound $A=\|\Lambda\|^{-1}\|g\|^2$.
\end{coro}
The statement is well-known, see \cite{CharlyBook} for the
elementary proof.
\begin{coro}
  Let $g_1,...,g_n,\gamma_1,...,\gamma_n\in S_0(\R^d)$,  then for
  the multi-window Gabor frame $S=\sum_{i=1}^nS_{g_i,\gamma_i}$
  the following are equivalent:
  \begin{enumerate}
    \item $S_{g_1,\gamma_1}+\cdots+S_{g_1,\gamma_1}=I$.
    \item $|\Lambda|^{-1}\big(\langle\gamma_1,\pi(\lambda^0)g_1\rangle+\cdots+
    \langle\gamma_1,\pi(\lambda^0)g_1\rangle\big)=\delta_{\lambda,0}$.
  \end{enumerate}
\end{coro}
The proof follows the same reasoning as for a single Gabor frame.
\section{Conclusions}

In the last decade operator algebra techniques have been of minor
interest in Gabor analysis. But in \cite{DLL95},\cite{Jan95} and
\cite{GL04} deep results about Gabor frames were obtained with the
help of operator algebras. We included our approach to the
Wexler-Raz biorthogonality principle as an indication for the
usefulness of Morita equivalence in Gabor analysis. In the
following we list some topics, where our approach gives new
insight, too.
\begin{enumerate}
  \item The original motivation for our study of Rieffel's results
        about Morita equivalence was the density result. There are different approaches
        to this important theorem \cite{DLL95},\cite{FK98} and \cite{Be04},
        which at the first sight seem unrelated. In \cite{Lu04a} we
        show that all these approaches cover different aspects of
        Morita equivalence between $C^*$-algebras generated by time-frequency shifts
        with respect to a lattice in the time-frequency plane.

  \item Our interpretation of Rieffel's construction of
        equivalence- bimodules for noncommutative tori in the notions of
        Gabor analysis enables us to answer the question posed by
        Manin on the connection between his quantum theta
        functions and the quantum theta vectors of Schwarz, see \cite{Lu04b}.
  \end{enumerate}
Other applications of Rieffel' setting yield new results on Feichtinger' conjecture 
and on the sturcture of multi-window Gabor frames, which is part of our current 
research.
\par
{\bf Acknowledgement:}
This investigations are part of the authors Ph.D. thesis
under the supervision of H.G. Feichtinger, whom I want to thank
for many helpful discussions and remarks on the content of the present paper. 
Finally, I want to express my gratitude to E. Matusiak and B. Scharinger for their 
careful reading of the manuscript.

\end{document}